\RequirePackage{fix-cm}
\documentclass[smallextended]{svjour3}       % onecolumn (second format)
\smartqed  % flush right qed marks, e.g. at end of proof
\usepackage{graphicx}

\usepackage{amsmath,amscd,amsfonts,amssymb,enumerate, mathrsfs}
\usepackage{color}
\usepackage[colorlinks]{hyperref}

\usepackage{floatrow}
\usepackage{multirow}
\usepackage[table]{xcolor}

\numberwithin{equation}{section}
\numberwithin{equation}{section}

\DeclareMathOperator{\RE}{Re}

\begin{document}
	
	\title{Radius of $\gamma$-Spirallikeness of order $\alpha$ for some Special functions
		%\thanks{The work of Kamajeet Gangania is supported by University Grant Commission, New-Delhi, India  under UGC-Ref. No.:1051/(CSIR-UGC NET JUNE 2017).}
	}
	%\subtitle{Do you have a subtitle?\\ If so, write it here}
	
	%\titlerunning{Short form of title}        % if too long for running head
	
	\author{ Sercan Kaz\i mo\u{g}lu       
		\and
		Kamaljeet Gangania$^*$
		 %etc.
	}
	
	\authorrunning{Sercan \and Kamaljeet} % if too long for running head
	
	\institute{	Sercan Kaz\i mo\u{g}lu \at
		\email{srcnkzmglu@gmail.com}  \\
		\emph{ Department of Mathematics, Faculty of Science and Literature, Kafkas University, Campus, 36100, Kars-Turkey}
		\and 
		Kamaljeet Gangania \at
		\email{gangania.m1991@gmail.com}             \\
		\emph{Department of Applied Mathematics, Delhi Technological University, Delhi--110042, India} %  if needed
		\and
		${}^*$ Corresponding author	    
	}
	
	\date{Received: date / Accepted: date}
	% The correct dates will be entered by the editor

	\maketitle

	\begin{abstract}
		In this paper, we establish the radius of $\gamma$-Spirallike of order $\alpha$ of certain well-known special functions. The main results of the paper are new and natural extensions of some known results.
		% \PACS{PACS code1 \and PACS code2 \and more}
	\end{abstract}
	
	\section{Introduction}
	Let $\mathcal{A}$ be the class of analytic functions normalized by the condition $f(0)=0=f'(0)-1$ in the unit disk $\mathbb{D}:=\mathbb{D}_{1}$, where   $\mathbb{D}_{r}:=\{z\in \mathbb{C}: |z|<r \}$.  We say that a function $f\in \mathcal{A}$ is $\gamma$-Spirallike of order $\alpha$ if and only if
	\begin{equation*}
	\RE\left(e^{-i\gamma} \frac{zf'(z)}{f(z)} \right) > \alpha \cos{\gamma},
	\end{equation*}
	where $\gamma \in \left(-\frac{\pi}{2},\frac{\pi}{2}\right)$ and $0\leq\alpha<1$. We denote the class of such functions by $\mathcal{S}_{p}^{\gamma}(\alpha)$. We also denote its convex analog, that is the class  $\mathcal{CS}_{p}^{\gamma}(\alpha)$ of convex $\gamma$-spirallike functions of order $\alpha$, which is defined below
	\begin{equation*}
	\RE\left(e^{-i\gamma} \left(1+\frac{zf''(z)}{f'(z)}\right) \right) > \alpha \cos{\gamma}.
	\end{equation*}
	
	The class $\mathcal{S}_{p}^{\gamma}(0)$ was introduced by Spacek~\cite{spacek-1933}. Each function in $\mathcal{CS}_{p}^{\gamma}(\alpha)$ is univalent in $\mathbb{D}$, but they do not necessarily be starlike. Further, it is worth to mention that for general values of $\gamma (|\gamma|<\pi/2)$, a function in $\mathcal{CS}_{p}^{\gamma}(0)$ need not be univalent in $\mathbb{D}$. For example: $f(z)=i(1-z)^i-i \in \mathcal{CS}_{p}^{\pi/4}(0)$, but not univalent. Indeed, $f \in \mathcal{CS}_{p}^{\gamma}(0)$ is univalent if $0<\cos\gamma<1/2$, see Robertson~\cite{Robertson-1969} and Pfaltzgraff~\cite{Pfaltzgraff-1975}.  Note that for $\gamma=0$, the classes $\mathcal{S}_{p}^{\gamma}(\alpha)$ and $\mathcal{CS}_{p}^{\gamma}(\alpha)$ reduce to the classes of starlike and convex functions of order $\alpha$, given by
	\begin{equation*}
	\RE\left(\frac{zf'(z)}{f(z)} \right) > \alpha \quad \text{and} \quad
	\RE\left(1+\frac{zf''(z)}{f'(z)}\right) > \alpha,
	\end{equation*}
	which we denote by $\mathcal{S}^*(\alpha)$ and $\mathcal{C}(\alpha)$, respectively. 
	
	In the recent past, connections between the special functions and their geometrical properties have been established in terms of radius problems \cite{abo-2018,bdoy-2016,Baric-2014,b-praj-2020,Baric-2015,btk-2018,Bricz-Rama,bulut-engel-2019,ErhanDenij2020,Kazimoglu-2022}. In this direction, behavior of the positive roots of a special function and the Laguerre-P\'{o}lya class play an evident role.  
	A real entire function $L$ maps real line into itself is said to be in the Laguerre-P\'{o}lya class $\mathcal{LP}$, if it can be expressed as follows:
	\begin{equation*}
	L(x)=cx^m e^{-ax^2+\beta x} \prod_{k\geq1}\left(1+\frac{x}{x_k} \right) e^{-\frac{x}{x_k}},
	\end{equation*}
	where $c,\beta,x_k\in \mathbb{R}$, $a\geq0$, $m\in \mathbb{N}\cup \{0 \}$ and $\sum {x_k}^{-2}< \infty$, see \cite{bdoy-2016}, \cite[p.~703]{lp}, \cite{Levin-1996} and the references therein. The class $\mathcal{LP}$ consists of entire functions which can be approximated by polynomials with only real zeros, uniformly on the compact sets of the complex plane and it is closed under differentiation.
	
The $\mathcal{S}^*(\alpha)$-radius, which is given below
$$\sup\{ r\in \mathbb{R}^+ : \RE\left(\frac{zg'(z)}{g(z)} \right) > \alpha, z\in \mathbb{D}_{r}  \}$$
and similarly, $\mathcal{C}(\alpha)$-radius has recently been obtained for some normalized forms of
 Bessel functions~\cite{abo-2018,Baric-2014,Baric-2015} (see Watson's treatise~\cite{watson-1944} for more on Bessel function), Struve functions~\cite{abo-2018,bdoy-2016}, Wright functions~\cite{btk-2018}, Lommel functions~\cite{abo-2018,bdoy-2016}, Legendre polynomials of odd degree~\cite{bulut-engel-2019} and Ramanujan type entire functions~\cite{ErhanDenij2020}. For their generalization to Ma-Minda classes~\cite{minda94} of starlike and convex functions, we refer to see~\cite{g-specialIJST,SG-2020}. 
 
 With the best of our knowledge, $\mathcal{S}_{p}^{\gamma}(\alpha)$-radius and $\mathcal{CS}_{p}^{\gamma}(\alpha)$-radius for special functions are not  handled till date. Therefore,
in this paper, we now aim to derive the radius of $\gamma$-Spirallike of order $\alpha$, which is given below
	\begin{equation*}
	R_{sp}(g) =\sup\left\{ r\in \mathbb{R}^+ : \RE\left(e^{-i\gamma} \frac{zg'(z)}{g(z)} \right) > \alpha \cos{\gamma}, z\in \mathbb{D}_{r} \right \}
	\end{equation*}
	and also the radius of convex $\gamma$-Spirallike of order $\alpha$, which is
	\begin{equation*}
	R_{sp}^{c}(g) = \sup\left\{ r\in \mathbb{R}^+ : \RE\left(e^{-i\gamma} \left(1+\frac{zg''(z)}{g'(z)}\right) \right) > \alpha \cos{\gamma}, z\in \mathbb{D}_{r} \right \}.
	\end{equation*}
	for the function $g$ in $\mathcal{A}$ to be a special function.

	\section{Wright functions}\label{sec-1}
	Let us consider the generalized Bessel function given by
	\begin{equation*}
	\Phi(\kappa, \delta, z) = \sum_{n\geq0} \frac{z^n}{n! \Gamma(n\kappa+ \delta)},
	\end{equation*}
	where $\kappa>-1$ and $z, \delta \in \mathbb{C}$, named after E. M. Wright. The function $\Phi$ is entire for $\kappa>-1$. From \cite[Lemma~1, p.~100]{btk-2018}, we have the Hadamard  factorization
	\begin{equation}\label{Had-wrt}
	\Gamma(\delta) \Phi(\kappa, \delta,-z^2)= \prod_{n\geq1}\left(1-\frac{z^2}{{\zeta}^{2}_{\kappa, \delta,n}}\right),
	\end{equation}
	where $\kappa, \delta>0$ and ${\zeta}_{\kappa, \delta,n}$ is the $n$-th positive root of $\Phi(\kappa, \delta,-z^2)$ and satisfies the interlacing property:
	\begin{equation}\label{w-roots}
	\breve{\zeta}_{\kappa, \delta, n}< {\zeta}_{\kappa, \delta, n} < \breve{\zeta}_{\kappa, \delta, n+1}< {\zeta}_{\kappa, \delta, n+1}, \quad (n\geq1)
	\end{equation}
	where $\breve{\zeta}_{\kappa, \delta,n}$ is the $n$-th positive root of the derivative of the function
	$$\Psi_{\kappa, \delta}(z)=z^{\delta}\Phi(\kappa, \delta,-z^2).$$ 
	Since $\Phi(\kappa, \delta, -z^2)\not \in \mathcal{A}$, therefore we choose the normalized Wright functions:
	\begin{equation}\label{w1}
	\left\{
	\begin{array}{lr}
	f_{\kappa, \delta}(z)=  \left[z^{\delta}\Gamma(\delta) \Phi(\kappa, \delta, -z^2)\right]^{1/\delta} \\
	
	g_{\kappa, \delta}(z)=   z\Gamma(\delta) \Phi(\kappa, \delta, -z^2) \\
	
	h_{\kappa, \delta}(z)=  z\Gamma(\delta) \Phi(\kappa, \delta, -z).
	\end{array}
	\right.
	\end{equation}
	
	For brevity, we write $W_{\kappa, \delta}(z):=\Phi(\kappa, \delta, -z^2)$. 
	\begin{theorem}\label{wright-phi}
		Let $\kappa, \delta>0$. The radius of $\gamma$-Spirallikeness for the functions $f_{\kappa, \delta}$, $g_{\kappa, \delta}$ and $h_{\kappa, \delta}$ are the smallest positive roots of the following equations:
		\begin{enumerate}
			\item [$(i)$]		$r W'_{\kappa, \delta}(r) +\delta\left(1-\alpha\right)\cos\gamma {W_{\kappa, \delta}(r)}=0$ 
			\item [$(ii)$]     	$r W'_{\kappa, \delta}(r) +\left(1-\alpha\right)\cos\gamma {W_{\kappa, \delta}(r)}=0$ 
			\item [$(iii)$]      $\sqrt{r} W'_{\kappa, \delta}(\sqrt{r}) +2\left(1-\alpha\right)\cos\gamma {W_{\kappa, \delta}(\sqrt{r})}=0$
		\end{enumerate}
		in $|z|< (0,{\zeta}_{\kappa, \delta,1})$, $(0,{\zeta}_{\kappa, \delta,1})$ and $(0,{\zeta}^2_{\kappa, \delta,1})$, respectively .
	\end{theorem}
	\begin{proof}
		Using \eqref{Had-wrt}, we obtain the following by the logarithmic differentiation of \eqref{w1}:
		\begin{equation}\label{w-sharp}
		\left\{
		\begin{array}{lr}
		\frac{zf'_{\kappa, \delta}(z)}{f_{\kappa, \delta}(z)} = 1+ \frac{1}{\delta} \frac{z W'_{\kappa, \delta}(z)}{W_{\kappa, \delta}(z)}=1-\frac{1}{\delta} \sum_{n\geq1}\frac{2z^2}{{\zeta}^2_{\kappa, \delta,n} -z^2} \\
		\frac{zg'_{\kappa, \delta}(z)}{g_{\kappa, \delta}(z)} = 1+ \frac{z W'_{\kappa, \delta}(z)}{W_{\kappa, \delta}(z)}=1- \sum_{n\geq1}\frac{2z^2}{{\zeta}^2_{\kappa, \delta,n} -z^2} \\
		\frac{zh'_{\kappa, \delta}(z)}{h_{\kappa, \delta}(z)} = 1+ \frac{1}{2} \frac{\sqrt{z} W'_{\kappa, \delta}(\sqrt{z})}{W_{\kappa, \delta}(\sqrt{z})}=1- \sum_{n\geq1}\frac{z}{{\zeta}^2_{\kappa, \delta,n} -z}. 
		\end{array}
		\right.
		\end{equation}
		We need to show that the following inequalities for 
		$\alpha \in \lbrack 0,1)$ and $\gamma \in \left(-\frac{\pi}{2},\frac{\pi}{2}\right),$
		\begin{equation}
		\operatorname{Re}\left( e^{-i\gamma}\frac{zf'_{\kappa, \delta}(z)}{f_{\kappa, \delta}(z)} \right) >\alpha\cos\gamma,
		~~~ \operatorname{Re}\left( e^{-i\gamma}\frac{zg'_{\kappa, \delta}(z)}{g_{\kappa, \delta}(z)} \right) >\alpha\cos\gamma
		\end{equation}
		and
		\begin{equation*}
		\operatorname{Re}\left(e^{-i\gamma} \frac{zh'_{\kappa, \delta}(z)}{h_{\kappa, \delta}(z)} \right) >\alpha\cos\gamma  \label{eq21}
		\end{equation*}
		are valid for $z\in {\mathbb{D}}_{r_{sp}(f_{\kappa, \delta})},~z\in {{\mathbb{D}}_{r_{sp}(g_{\kappa, \delta})}}$ and $z\in {{\mathbb{D}}_{r_{sp}(h_{\kappa, \delta})}}$ respectively, and each of the above inequalities does not hold
		in larger disks.
		It is known \cite{Deniz-2017} that if $z\in \mathbb{C}$ and $\lambda \in \mathbb{R}$ are such that $ \left\vert z\right\vert\leq r<\lambda ,$ then
		\begin{equation}
		\operatorname{Re}\left( \frac{z}{\lambda -z}\right) \leq \left\vert \frac{z}{\lambda -z}\right\vert  \leq \frac{|z|}{\lambda -\left\vert z\right\vert }.  \label{eq23}
		\end{equation}
		Then the inequality
		\begin{equation*}
		\operatorname{Re}\left( \frac{z^{2}}{ {\zeta}^2_{\kappa, \delta,n}-z^{2}}\right) \leq \left\vert \frac{z^{2}}{ {\zeta}^2_{\kappa, \delta,n}-z^{2}}\right\vert \leq \frac{\left\vert z\right\vert ^{2}}{ {\zeta}^2_{\kappa, \delta,n}-\left\vert z\right\vert ^{2}}
		\end{equation*}%
		holds for every $\left|z\right|<{\zeta}_{\kappa, \delta,1}. $ Therefore, from (\ref{w-sharp}) and (\ref{eq23}), we have 
		\begin{eqnarray}
		\operatorname{Re}\left(  e^{-i\gamma}\frac{zf'_{\kappa, \delta}(z)}{f_{\kappa, \delta}(z)}\right) 
		&=& \operatorname{Re}\left( e^{-i\gamma}\right)-\frac{1}{\delta}\operatorname{Re}\left(e^{-i\gamma}\sum\limits_{n\geq 1} \frac{2z^{2}}{ {\zeta}^2_{\kappa, \delta,n}-z^{2}}\right)  \nonumber\\
		&\geq& \cos\gamma- \frac{1}{\delta}\left|e^{-i\gamma}\sum\limits_{n\geq 1} \frac{2z^{2}}{ {\zeta}^2_{\kappa, \delta,n}-z^{2}}\right|
		\geq  \cos\gamma-\frac{1}{\delta }\sum\limits_{n\geq 1}\frac{2|z|^{2}}{ {\zeta}^2_{\kappa, \delta,n}- |z|^{2}}   \nonumber \\
		&=&\frac{|z| f_{\kappa, \delta}(|z| )}{f_{\kappa, \delta}(|z|)}+\cos\gamma-1. \label{w-mod}
		\end{eqnarray}
		Equality in the each of the above inequalities \eqref{w-mod} holds when $z=r$. Thus, for $r\in \left(0,{\zeta}_{\kappa, \delta,1} \right)$ it follows that 
		\begin{equation*}
		\inf_{z\in \mathbb{D}_r } \left\lbrace \operatorname{Re} \left( e^{-i\gamma}\frac{zf'_{\kappa, \delta}(z)}{f_{\kappa, \delta}(z)}-\alpha\cos\gamma \right) \right\rbrace= \frac{\left\vert z\right\vert f'_{\kappa, \delta}(\left\vert	z\right\vert )}{f_{\kappa, \delta}(\left\vert z\right\vert )}+\left(1-\alpha\right)\cos\gamma-1.
		\end{equation*}
		Now, the mapping $\Theta :\left(0,{\zeta}_{\kappa, \delta,1} \right)\longrightarrow \mathbb{R}$ defined by 
		\begin{equation*}
		\Theta(r)= \frac{rf'_{\kappa, \delta}(r)}{f_{\kappa, \delta}(r)} +\left(1-\alpha\right)\cos\gamma-1=\left(1-\alpha\right)\cos\gamma-\frac{1}{\delta}
		\sum\limits_{n\geq 1}\left( \frac{2r^{2}}{ {\zeta}_{\kappa, \delta,n}^{2}-r^{2}}\right).
		\end{equation*}
		is strictly decreasing since
		\begin{equation*}
		\Theta^\prime(r)= -\frac{1}{\delta}\sum\limits_{n\geq 1}\left( \frac{4r{\zeta}_{\kappa, \delta,n}}{ \left({\zeta}_{\kappa, \delta,n}^{2}-r^{2}\right)^2}\right)<0
		\end{equation*}
		for all $\delta>0.$ On the other hand, since 
		$$\lim_{r\searrow0}\Theta(r)=\left(1-\alpha\right)\cos\gamma>0  \text{ \  and \ } \lim_{r\nearrow {\zeta}_{\kappa, \delta,1}}\Theta(r)=-\infty,$$
		in view of the minimum principle for harmonic functions imply that the corresponding inequality for $f_{\kappa, \delta}$ in (\ref{eq21})
		for $\delta>0$ holds if and only if $z\in \mathbb{D}_{r_{sp}(f_{\kappa, \delta})},$ where ${r_{sp}(f_{\kappa, \delta})}$ is the smallest positive root of equation 
		\begin{equation*}
		\frac{rf'_{\kappa, \delta}(r)}{f_{\kappa, \delta}(r)}=1-\left(1-\alpha\right)\cos\gamma
		\end{equation*}
		which is equivalent to 
		\begin{equation*}
		\frac{1}{\delta} \frac{z W'_{\kappa, \delta}(z)}{W_{\kappa, \delta}(z)}=-\left(1-\alpha\right)\cos\gamma, 
		\end{equation*}
		situated in $\left(0,{\zeta}_{\kappa, \delta,1} \right).$ Reasoning along the same lines, proofs of the other parts follows. \qed
	\end{proof}
	\begin{remark}
		Taking $\gamma=0$ in Theorem~\ref{wright-phi} yields \cite[Theorem~1]{btk-2018}. 
	\end{remark}
	
	In the following, we deal with convex analogue of the class of $\gamma$-spiralllike functions of order $\alpha$.
	\begin{theorem}\label{wright-conx}
		Let $\kappa, \delta>0$ and the functions $f_{\kappa, \delta}$, $g_{\kappa, \delta}$ and $h_{\kappa, \delta}$ as given in \eqref{w1}. Then
		\begin{enumerate}
			\item [(i)] the radius $R_{sp}^{c}(f_{\kappa, \delta})$ is the smallest positive root of the equation 
			$$\frac{r {\Psi}^{''}_{\kappa, \delta}(r)}{{\Psi}^{'}_{\kappa, \delta}(r)}+ \left(\frac{1}{\delta}-1\right)\frac{r {\Psi}^{'}_{\kappa, \delta}(r)}{{\Psi}_{\kappa, \delta}(r)} +(1-\alpha)\cos\gamma=0.$$
			
			\item [(ii)] the radius $R_{sp}^{c}(g_{\kappa, \delta})$ is the smallest positive root of the equation 
			$$rg''_{\kappa, \delta}(r) +(1-\alpha)\cos\gamma g'_{\kappa, \delta}(r)=0.$$
			
			\item [(iii)] the radius $R_{sp}^{c}(h_{\kappa, \delta})$ is the smallest positive root of the equation 
			$$rh''_{\kappa, \delta}(r) +(1-\alpha)\cos\gamma h'_{\kappa, \delta}(r)=0.$$
		\end{enumerate}
	\end{theorem}
	\begin{proof}
		We first prove the part $(i)$. From \eqref{Had-wrt}, \eqref{w1} and using the Hadamard representation $\Gamma(\delta) {\Psi}^{'}_{\kappa, \delta}(z)=\delta z^{\delta-1} \prod_{n\geq1}\left(1-\frac{z^2}{ \breve{\zeta}^2_{\kappa, \delta, n}} \right)$, (see \cite[Eq.~7]{btk-2018}), we have
		\begin{align*}
		1+\frac{z f''_{\kappa, \delta}(z)}{f'_{\kappa, \delta}(z)}&=1+\frac{z {\Psi}^{''}_{\kappa, \delta}(z)}{{\Psi}^{'}_{\kappa, \delta}(z)}+ \left(\frac{1}{\delta}-1\right)\frac{z {\Psi}^{'}_{\kappa, \delta}(z)}{{\Psi}_{\kappa, \delta}(z)}\\
		&= 1-\sum_{n\geq1}\frac{2z^2}{\breve{\zeta}^2_{\kappa, \delta, n} -z^2} -\left(\frac{1}{\delta}-1\right)\sum_{n\geq1}\frac{2z^2}{	{\zeta}^2_{\kappa, \delta, n}-z^2}
		\end{align*}
		and for $\delta>1$, using the following inequality of \cite{Deniz-2017}: 
		\begin{equation}\label{firstnorm}
		\left|\frac{z}{y-z}-\lambda \frac{z}{x-z}\right| \leq \frac{|z|}{y-|z|}-\lambda \frac{|z|}{x-|z|},\quad (x>y>r\geq|z|)
		\end{equation}
		with $\lambda=1-1/\delta$, we get
		\begin{align*}
		\left|\frac{z f''_{\kappa, \delta}(z)}{f'_{\kappa, \delta}(z)} \right|&\leq
		-\frac{r f''_{\kappa, \delta}(r)}{f'_{\kappa, \delta}(r)} =-\frac{r {\Psi}^{''}_{\kappa, \delta}(r)}{{\Psi}^{'}_{\kappa, \delta}(r)}- \left(\frac{1}{\delta}-1\right)\frac{r {\Psi}^{'}_{\kappa, \delta}(r)}{{\Psi}_{\kappa, \delta}(r)}.
		\end{align*}
		Also, using the inequality $||x|-|y||\leq |x-y|$ and the relation in \eqref{w-roots}, we see that for $\delta>0$
		\begin{equation*}
		\left|\frac{z f''_{\kappa, \delta}(z)}{f'_{\kappa, \delta}(z)}\right| \leq -\frac{r f''_{\kappa, \delta}(r)}{f'_{\kappa, \delta}(r)},
		\end{equation*}
		holds in$|z|=r< \breve{\zeta}_{\kappa, \delta,1}$. Therefore, we have 
		\begin{align}
		& \operatorname{Re}\left(  e^{-i\gamma}\left(1+\frac{zf''_{\kappa, \delta}(z)}{f'_{\kappa, \delta}(z)}\right)\right) \nonumber\\
		&= \operatorname{Re}(e^{-i\gamma}) -\operatorname{Re}\left(e^{-i\gamma}\left(\sum_{n\geq1}\frac{2z^2}{\breve{\zeta}^2_{\kappa, \delta, n} -z^2} +\left(\frac{1}{\delta}-1\right)\sum_{n\geq1}\frac{2z^2}{	{\zeta}^2_{\kappa, \delta, n}-z^2} \right) \right)  \nonumber\\
		&\geq \cos\gamma- \left|\sum_{n\geq1}\frac{2z^2}{\breve{\zeta}^2_{\kappa, \delta, n} -z^2} +\left(\frac{1}{\delta}-1\right)\sum_{n\geq1}\frac{2z^2}{	{\zeta}^2_{\kappa, \delta, n}-z^2} \right| \nonumber\\
		&\geq  \cos\gamma+\frac{r f''_{\kappa, \delta}(r)}{f'_{\kappa, \delta}(r)} \label{w-mod}
		\end{align}
		hold for $\delta>1$. Observe that these inequalities also hold for $\delta>0$. Equality in the each of the above inequalities \eqref{w-mod} holds when $z=r$. Thus, for $r\in ( 0,\breve{\zeta}_{\kappa, \delta,1} )$ it follows that 
		\begin{equation*}
		\inf_{z\in \mathbb{D}_r } \left\lbrace \operatorname{Re} \left( e^{-i\gamma}\left(1+\frac{zf''_{\kappa, \delta}(z)}{f'_{\kappa, \delta}(z)}\right)-\alpha\cos\gamma \right) \right\rbrace= \left(1-\alpha\right)\cos\gamma+\frac{\left\vert z\right\vert f''_{\kappa, \delta}(\left\vert	z\right\vert )}{f'_{\kappa, \delta}(\left\vert z\right\vert )}.
		\end{equation*}
		Now, the proof of part $(i)$ follows on similar lines as of Theorem~\ref{wright-phi}.
		
		For the other parts, note that the functions $g_{\kappa, \delta}$ and $h_{\kappa, \delta}$ belong to the Laguerre-P\'{o}lya class $\mathcal{LP}$, which is closed under differentiation, their derivatives $g'_{\kappa, \delta}$ and $h'_{\kappa, \delta}$ also belong to $\mathcal{LP}$ and the zeros are real. Thus assuming $\tau_{\kappa, \delta,n}$  and $\eta_{\kappa, \delta,n}$ are the positive zeros of $g'_{\kappa, \delta}$ and $h'_{\kappa, \delta}$, respectively, we have the following representations:
		\begin{align*}
		g'_{\kappa, \delta}(z)=\prod_{n\geq1}\left(1-\frac{z^2}{{\tau}^2_{\kappa, \delta,n}}\right) \quad
		\text{and} \quad
		h'_{\kappa, \delta}(z)=\prod_{n\geq1}\left(1-\frac{z}{\eta_{\kappa, \delta,n}}\right),
		\end{align*}
		which yield
		\begin{align*}
		1+\frac{zg''_{\kappa, \delta}(z)}{g'_{\kappa, \delta}(z)}=1-\sum_{n\geq1}\frac{2z^2}{{\tau}^2_{\kappa, \delta,n}-z^2} \quad
		\text{and}\quad
		1+\frac{zh''_{\kappa, \delta}(z)}{h'_{\kappa, \delta}(z)}=1-\sum_{n\geq1}\frac{z}{{\eta}_{\kappa, \delta,n}-z}.
		\end{align*}
		Further, reasoning along the same lines as in Theorem~\ref{wright-phi}, the result follows at once. \qed
	\end{proof}
	\begin{remark}
		 Taking $\gamma=0$ in Theorem~\ref{wright-conx} yields \cite[Theorem~5]{btk-2018}. 
	\end{remark}
	
	\section{Mittag-Leffler functions}\label{sec-2}
	In 1971, Prabhakar~\cite{prabha-1971} introduced the following function
	\begin{equation*}
	M(\mu,\nu, a, z):= \sum_{n\geq0} \frac{(a)_{n} z^n}{n! \Gamma(\mu n+\nu)},
	\end{equation*}
	where $(a)_{n}=\Gamma(a+n)/\Gamma(a)$ denotes the Pochhammer symbol and $\mu, \nu,a>0$. The functions $M(\mu,\nu, 1, z)$ and $M(\mu, 1, 1, z)$ were introduced and studied by Wiman and Mittag-Leffler, respectively. Now let us consider the set $ W_{b}= A(W_{c})\cup B(W_{c})$, where
	\begin{equation*}
	W_{c}:= \left\{  \left(\frac{1}{\mu},\nu\right): 1<\mu<2, \nu\in [\mu-1,1]\cup[\mu,2]   \right\}
	\end{equation*}
	and denote by $W_{i}$, the smallest set containing $W_{b}$ and invariant under the transformations $A$, $B$ and $C$ mapping the set
	$\{(\tfrac{1}{\mu},\nu): \mu>1, \nu>0\}$ into itself and are defined as:
	\begin{align*}
	&A: (\tfrac{1}{\mu},\nu)\rightarrow (\tfrac{1}{2\mu},\nu),\quad B: (\tfrac{1}{\mu},\nu)\rightarrow (\tfrac{1}{2\mu},\mu+\nu), \\
	&C: (\tfrac{1}{\mu},\nu)\rightarrow 
	\left\{
	\begin{array}
	{lr}
	(\tfrac{1}{\mu},\nu-1),     & \text{if}\; \nu>1; \\
	(\tfrac{1}{\mu},\nu),     & \text{if}\; 0<\nu\leq1.
	\end{array}
	\right.
	\end{align*}
	Kumar and Pathan \cite{pathan-2016} proved that if $(\tfrac{1}{\mu},\nu)\in W_{i}$ and $a>0$, then all zeros of $M(\mu,\nu, a, z)$ are real and negative. From \cite[Lemma~1, p.~121]{b-praj-2020}, we see that if  $(\tfrac{1}{\mu},\nu)\in W_{i}$ and $a>0$, then the function $M(\mu,\nu, a, -z^2)$ has infinitley many zeros, which are all real and have the following representation:
	\begin{equation*}
	\Gamma(\nu)M(\mu,\nu, a, -z^2)=\prod_{n\geq1}\left(1-\frac{z^2}{{\lambda}^2_{\mu,\nu,a,n}}\right),
	\end{equation*}
	where ${\lambda}_{\mu,\nu,a,n}$ is the $n$-th positive zero of $M(\mu,\nu, a, -z^2)$ and satisfy the interlacing relation
	\begin{equation*}
	{\xi}_{\mu,\nu, a,n}< {\lambda}_{\mu,\nu,a,n}< {\xi}_{\mu,\nu, a,n+1} < {\lambda}_{\mu,\nu,a,n+1} \quad (n\geq1),
	\end{equation*}
	where ${\xi}_{\mu,\nu, a,n}$ is the $n$-th positive zero of the derivative of $ z^{\nu}M(\mu,\nu, a, -z^2)$. Since $M(\mu,\nu, a, -z^2)\not\in \mathcal{A}$, therefore we consider the following normalized forms (belong to the Laguerre-P\'{o}lya class):
	\begin{equation}\label{mittag1}
	\left\{
	\begin{array}{lr}
	f_{\mu, \nu,a}(z) = \left[z^{\nu}\Gamma(\nu) M(\mu, \nu, a, -z^2)\right]^{1/\nu}, \\
	g_{\mu, \nu,a}(z) = z\Gamma(\nu) M(\mu, \nu, a, -z^2) \\
	h_{\mu, \nu,a}(z) = z\Gamma(\nu) M(\mu, \nu,a, -z).
	\end{array}
	\right.
	\end{equation}
	
	For brevity, write $L_{\mu,\nu,a}(z):=M(\mu,\nu, a, -z^2)$. Now proceeding similarly as in Section~\ref{sec-1}, we obtain the following results:
	\begin{theorem}\label{mittag-phi}
		Let $(\tfrac{1}{\mu},\nu)\in W_{i}$, $a>0$. Then the radius of $\gamma$-Spirallikeness of order $\alpha$
		for the functions $f_{\mu,\nu, a}$, $g_{\mu,\nu, a}$ and $h_{\mu,\nu, a}$ given by \eqref{mittag1} are the smallest positive roots of the following equations:
		\begin{enumerate}
			\item [$(i)$]		$r L'_{\mu, \nu,a}(r) +\delta\left(1-\alpha\right)\cos\gamma {L_{\mu, \nu,a}(r)}=0$ 
			\item [$(ii)$]     	$r L'_{\mu, \nu,a}(r) +\left(1-\alpha\right)\cos\gamma {L_{\mu, \nu,a}(r)}=0$ 
			\item [$(iii)$]      $\sqrt{r} L'_{\mu, \nu,a}(\sqrt{r}) +2\left(1-\alpha\right)\cos\gamma {L'_{\mu, \nu,a}(\sqrt{r})}=0$
		\end{enumerate}
		in $|z|< (0,{\lambda}_{\mu, \nu,a,1})$, $(0,{\lambda}_{\mu, \nu,a,1})$ and $(0,{\lambda}^2_{\mu, \nu,a,1})$, respectively .		
	\end{theorem}
	\begin{proof} 
		Using \eqref{mittag1}, we obtain after the logarithmic differentiation:
		\begin{equation}\label{mittag-starexpress}
		\left\{
		\begin{array}{lr}
		\frac{zf'_{\mu,\nu, a}(z)}{f_{\mu,\nu, a}(z)} =1-\frac{1}{\nu} \sum_{n\geq1}\frac{2z^2}{{\lambda}^2_{\mu,\nu,a,n} -z^2} \\
		\frac{zg'_{\mu,\nu, a}(z)}{g_{\mu,\nu, a}(z)} =1- \sum_{n\geq1}\frac{2z^2}{{\lambda}^2_{\mu,\nu,a,n} -z^2} \\
		\frac{zh'_{\mu,\nu, a}(z)}{h_{\mu,\nu, a}(z)} =1- \sum_{n\geq1}\frac{z}{{\lambda}^2_{\mu,\nu,a,n} -z}. 
		\end{array}
		\right.
		\end{equation}
		We need to show that the following inequalities for 
		$\alpha \in \lbrack 0,1)$ and $\gamma \in \left(-\frac{\pi}{2},\frac{\pi}{2}\right),$
		\begin{equation}\label{eq21-mittag}
		\operatorname{Re}\left( e^{-i\gamma}\frac{zf'_{\mu,\nu, a}(z)}{f_{\mu,\nu, a}(z)} \right) >\alpha\cos\gamma,
		~~~ \operatorname{Re}\left( e^{-i\gamma}\frac{zg'_{\mu,\nu, a}(z)}{g_{\mu,\nu, a}(z)} \right) >\alpha\cos\gamma
		\end{equation}
		and
		\begin{equation*}
		\operatorname{Re}\left(e^{-i\gamma} \frac{zh'_{\mu,\nu, a}(z)}{h_{\mu,\nu, a}(z)} \right) >\alpha\cos\gamma 
		\end{equation*}
		are valid for $z\in {\mathbb{D}}_{r_{sp}(f_{\mu, \nu,a})},~z\in {{\mathbb{D}}_{r_{sp}(g_{\mu, \nu,a})}}$ and $z\in {{\mathbb{D}}_{r_{sp}(h_{\mu, \nu,a})}}$ respectively, and each of the above inequalities does not hold in larger disks.  Since using~\eqref{eq23}
		\begin{equation}\label{eq22-mittag}
		\operatorname{Re}\left( \frac{z^{2}}{ {\lambda}^2_{\mu,\nu,a,n}-z^{2}}\right) \leq \left\vert \frac{z^{2}}{ {\lambda}^2_{\mu,\nu,a,n}-z^{2}}\right\vert \leq \frac{\left\vert z\right\vert ^{2}}{ {\lambda}^2_{\mu,\nu,a,n}-\left\vert z\right\vert ^{2}}
		\end{equation}%
		holds for every $\left|z\right|<{\lambda}_{\mu,\nu,a,1}. $ Therefore, from (\ref{mittag-starexpress}) and (\ref{eq22-mittag}), we have 
		\begin{eqnarray}
		\operatorname{Re}\left( e^{-i\gamma}\frac{zf'_{\mu,\nu, a}(z)}{f_{\mu,\nu, a}(z)} \right)
		&=& \operatorname{Re}\left( e^{-i\gamma}\right)-\frac{1}{\nu}\operatorname{Re}\left(e^{-i\gamma}\sum_{n\geq1}\frac{2z^2}{{\lambda}^2_{\mu,\nu,a,n} -z^2}\right)  \nonumber\\
		&\geq& \cos\gamma- \frac{1}{\nu}\left|e^{-i\gamma}\sum_{n\geq1}\frac{2z^2}{{\lambda}^2_{\mu,\nu,a,n} -z^2}\right| \nonumber \\
		&\geq&  \cos\gamma-\frac{1}{\nu }\sum\limits_{n\geq 1}\frac{2|z|^{2}}{ {\lambda}^2_{\mu,\nu,a,n}- |z|^{2}}   \nonumber \\
		&=&\frac{|z| f'_{\mu,\nu, a}(|z| )}{f_{\mu,\nu, a}(|z|)}+\cos\gamma-1. \label{mittag-mod}
		\end{eqnarray}
		Equality in the each of the above inequalities \eqref{mittag-mod} holds when $z=r$. Thus, for $r\in \left(0,{\lambda}_{\mu,\nu,a,1} \right)$ it follows that 
		\begin{equation*}
		\inf_{z\in \mathbb{D}_r } \left\lbrace \operatorname{Re} \left( e^{-i\gamma}\frac{zf'_{\mu,\nu, a}(z)}{f_{\mu,\nu, a}(z)}-\alpha\cos\gamma \right) \right\rbrace= \frac{|z| f'_{\mu,\nu, a}(|z| )}{f_{\mu,\nu, a}(|z|)}+\left(1-\alpha\right)\cos\gamma-1.
		\end{equation*}
		Now, the mapping $\Theta :\left(0,{\lambda}_{\mu, \nu,a,1} \right)\longrightarrow \mathbb{R}$ defined by 
		\begin{equation*}
		\Theta(r)= \frac{rf'_{\mu, \nu, a}(r)}{f_{\mu, \nu, a}(r)} +\left(1-\alpha\right)\cos\gamma-1=\left(1-\alpha\right)\cos\gamma-\frac{1}{\nu}
		\sum\limits_{n\geq 1}\left( \frac{2r^{2}}{ {\lambda}_{\mu, \nu,a,n}^{2}-r^{2}}\right).
		\end{equation*}
		is strictly decreasing since
		\begin{equation*}
		\Theta^\prime(r)= -\frac{1}{\nu}\sum\limits_{n\geq 1}\left( \frac{4r{\lambda}_{\mu, \nu,a,n}}{ \left({\lambda}_{\mu, \nu,a,n}^{2}-r^{2}\right)^2}\right)<0
		\end{equation*}
		for all $\nu>0.$ On the other hand, since 
		$$\lim_{r\searrow0}\Theta(r)=\left(1-\alpha\right)\cos\gamma>0  \text{ \  and \ } \lim_{r\nearrow {\lambda}_{\mu, \nu,a,1}}\Theta(r)=-\infty,$$
		in view of the minimum principle for harmonic functions imply that the corresponding inequality for $f_{\mu, \nu,a}$ in (\ref{eq21-mittag})
		for $\nu>0$ holds if and only if $z\in \mathbb{D}_{r_{sp}(f_{\mu, \nu,a})},$ where ${r_{sp}(f_{\mu, \nu,a})}$ is the smallest positive root of equation 
		\begin{equation*}
		\frac{rf'_{\mu, \nu,a}(r)}{f_{\mu, \nu,a}(r)}=1-\left(1-\alpha\right)\cos\gamma
		\end{equation*}
		which is equivalent to 
		\begin{equation*}
		\frac{1}{\nu} \frac{z L'_{\mu, \nu,a}(z)}{L_{\mu, \nu,a}(z)}=-\left(1-\alpha\right)\cos\gamma, 
		\end{equation*}
		situated in $\left(0,{\lambda}_{\mu, \nu,1} \right).$ Reasoning along the same lines, proofs of the other parts follows. \qed
	\end{proof}
	\begin{remark}
			Taking $\gamma=0$ in Theorem~\ref{mittag-phi} yields \cite[Theorem~1]{b-praj-2020}.
	\end{remark}
	
	In the following, we derive the result for the convex analog proceedings on similar lines as Theorem~\ref{wright-conx}.
	\begin{theorem}\label{mittag-c}
		Let $(\tfrac{1}{\mu},\nu)\in W_{i}$, $a>0$. Let the functions $f_{\mu,\nu, a}$, $g_{\mu,\nu, a}$ and $h_{\mu,\nu, a}$  be given by \eqref{mittag1}.Then
		\begin{enumerate}
			\item [(i)] the radius $R_{sp}^{c}(f_{\mu,\nu, a})$ is the smallest positive root of the equation 
			$$rf''_{\mu,\nu, a}(r) +(1-\alpha)\cos\gamma f'_{\mu,\nu, a}(r)=0.$$
			
			\item [(ii)] the radius $R_{sp}^{c}(g_{\mu,\nu, a})$ is the smallest positive root of the equation 
			$$rg''_{\mu,\nu, a}(r) +(1-\alpha)\cos\gamma g'_{\mu,\nu, a}(r)=0.$$
			
			\item [(iii)] the radius $R_{sp}^{c}(h_{\mu,\nu, a})$ is the smallest positive root of the equation 
			$$rh''_{\mu,\nu, a}(r) -(1-\alpha)\cos\gamma h'_{\mu,\nu, a}(r)=0.$$
		\end{enumerate}
	\end{theorem}
	\begin{remark}
			Taking $\gamma=0$ in Theorem~\ref{mittag-c} yields \cite[Theorem~3]{b-praj-2020}.
	\end{remark}	
	
	\section{Legendre polynomials}
	The Legendre polynomials $P_{n}$ are the solutions of the Legendre differential equation
	$$((1-z^2)P'_{n}(z))'+n(n+1)P_{n}(z)=0,$$
	where $n\in \mathbb{Z}^{+}$ and using Rodrigues formula, $P_{n}$ can be represented in the form:
	$$P_{n}(z)=\frac{1}{2^n n!}\frac{d^n(z^2-1)^n}{dz^n}$$
	and it also satisfies the geometric condition $P_n(-z)=(-1)^n P_{n}(z)$. Moreover, the odd degree Legendre polynomials $P_{2n-1}(z)$  have only real roots which satisfy 
	\begin{equation}\label{legdroot}
	0=z_0<z_1<\cdots<z_{n-1}\quad\text{or}\quad -z_1>\cdots>-z_{n-1}.
	\end{equation}
	Thus the normalized form is as follows:
	\begin{equation}\label{legd1}
	\mathcal{P}_{2n-1}(z):=\frac{P_{2n-1}(z)}{P'_{2n-1}(0)}=z+\sum_{k=2}^{2n-1}a_{k}z^{k}=a_{2n-1}z\prod_{k=1}^{n-1}(z^2-z^2_{k}).
	\end{equation}

	\begin{theorem}\label{Leg-c}
		Let $\mathcal{P}_{2n-1}$ be given by \eqref{legd1}. Then
		\begin{enumerate}
			\item [$(i)$] the radius $R_{sp}^{c}(\mathcal{P}_{2n-1})$ is the smallest positive root of the equation 
			$${r \mathcal{P}''_{2n-1}(r)} +(1-\alpha) {\mathcal{P}'_{2n-1}(r)}=0.$$
			
			\item [$(ii)$] the radius of $\gamma$-Spirallikeness of order $\alpha$ for the normalized Legendre polynomial of odd degree is given by the smallest positive root of the equation
			\begin{equation*}
			{r \mathcal{P}'_{2n-1}(r)} +(1-\alpha) {\mathcal{P}_{2n-1}(r)}=0.
			\end{equation*} 
		\end{enumerate} 
	\end{theorem}
	\begin{proof}
		We prove first part and second part follows on same lines. From \eqref{legd1}, upon the logarithmic differentiation, we have
		\begin{align*}
		1+\frac{z \mathcal{P}''_{2n-1}(z)}{\mathcal{P}'_{2n-1}(z)}=\frac{z \mathcal{P}'_{2n-1}(z)}{\mathcal{P}_{2n-1}(z)} -\frac{\sum_{k=1}^{n-1}\frac{4z^2_k z^2}{(z^2_k -z^2)^2}}{\frac{z \mathcal{P}'_{2n-1}(z)}{\mathcal{P}_{2n-1}(z)}},
		\end{align*}
		where
		\begin{align*}
		\frac{z \mathcal{P}'_{2n-1}(z)}{\mathcal{P}_{2n-1}(z)}=1-\sum_{k=1}^{n-1}\frac{2z^2}{z^2_k-z^2}.
		\end{align*}
		Further, after using the inequality $||x|-|y||\leq |x-y|$ and \eqref{legdroot} for $ |z|=r<z_1$, we see that
		\begin{align*}
		&\operatorname{Re}\left( e^{-i\gamma}\left(1+\frac{z \mathcal{P}''_{2n-1}(z)}{\mathcal{P}'_{2n-1}(z)}\right) \right)\\
		=& \operatorname{Re}\left( e^{-i\gamma}\right)-\operatorname{Re}\left(e^{-i\gamma}\sum_{k=1}^{n-1}\frac{2z^2}{z^2_k-z^2}\right)
		-\operatorname{Re}\left(e^{-i\gamma}\frac{\sum_{k=1}^{n-1}\frac{4z^2_k z^2}{(z^2_k -z^2)^2}}{1-\sum_{k=1}^{n-1}\frac{2z^2}{z^2_k-z^2}}\right)  \nonumber\\
		\geq& \cos\gamma-\left|e^{-i\gamma}\sum_{k=1}^{n-1}\frac{2z^2}{z^2_k-z^2}\right|
		-\left|e^{-i\gamma}\frac{\sum_{k=1}^{n-1}\frac{4z^2_k z^2}{(z^2_k -z^2)^2}}{1-\sum_{k=1}^{n-1}\frac{2z^2}{z^2_k-z^2}}\right|   \nonumber \\
		=& \cos\gamma-\sum_{k=1}^{n-1}\frac{2r^2}{z^2_k-r^2}-\frac{\sum_{k=1}^{n-1}\frac{4z^2_k r^2}{(z^2_k -r^2)^2}}{1-\sum_{k=1}^{n-1}\frac{2r^2}{z^2_k-r^2}} 
		=\cos\gamma + \frac{r \mathcal{P}''_{2n-1}(r)}{\mathcal{P}'_{2n-1}(r)}.
		\end{align*}
		Further, with similar reasoning as Theorem~\ref{mittag-phi}, result follows. \qed 
	\end{proof}
	\begin{remark}
		Taking $\gamma=0$ in Theorem~\ref{Leg-c} yields \cite[Theorem~2.2]{bulut-engel-2019} and \cite[Theorem~2.1]{bulut-engel-2019}.
	\end{remark}	
	
	\section{ Lommel functions}
	The Lommel function $\mathcal{L}_{u,v}$ of first kind is a particular solution of the second-order inhomogeneous Bessel differential equation $$z^2w''(z)+zw'(z)+(z^2-{v}^2)w(z)=z^{u+1},$$
	where $u\pm v\notin \mathbb{Z}^{-}$ and is given by
	$$\mathcal{L}_{u,v}=\frac{z^{u+1}}{(u-v+1)(u+v+1)}{}_1F_2\left(1;\frac{u-v+3}{2},\frac{u+v+3}{2};-\frac{z^2}{4}\right),$$
	where $\frac{1}{2}(-u\pm v-3)\notin \mathbb{N}$ and ${}_1 F_{2}$ is a hypergeometric function. Since it is not normalized, therefore we consider the following three normalized functions involving $\mathcal{L}_{u,v}$ :
	\begin{equation}\label{fL}
	\left\{
	\begin{array}{lr}
	f_{u,v}(z)=((u-v+1)(u+v+1)\mathcal{L}_{u,v}(z))^{\tfrac{1}{u+1}},\\
	g_{u,v}(z)=(u-v+1)(u+v+1)z^{-u}\mathcal{L}_{u,v}(z),\\
	h_{u,v}(z)=(u-v+1)(u+v+1)z^{(1-u)/2}\mathcal{L}_{u,v}(\sqrt{z}).
	\end{array}
	\right.
	\end{equation}
	Authors in \cite{abo-2018,bdoy-2016} and \cite{Bricz-Rama} proved the radius of starlikeness and convexity for the following normalized functions expressed in terms of $\mathcal{L}_{u-\tfrac{1}{2},\tfrac{1}{2}}$:
	\begin{equation}\label{lomel-normalized}
	f_{u-\tfrac{1}{2},\tfrac{1}{2}}(z),\quad g_{u-\tfrac{1}{2},\tfrac{1}{2}}(z) \quad\text{and}\quad h_{u-\tfrac{1}{2},\tfrac{1}{2}}(z),
	\end{equation} 
	where $0\neq u\in (-1,1)$.
	
	For brevity, we write these as $f_{u}, g_{u}$ and $h_{u}$, respectively and $\mathcal{L}_{u-\tfrac{1}{2},\tfrac{1}{2}}=\mathcal{L}_{u}$.
	
	\begin{theorem}\label{lommeltheorem}
		Let $u\in (-1,1)$, $u\neq0$. Let the functions $f_{u}, g_{u}$ and $h_{u}$ be given by \eqref{lomel-normalized}. Then
		\begin{enumerate}
			\item [(i)] the radius $R_{sp}^{c}(f_{u})$ is the smallest positive root of the equation 
			$$rf''_{u}(r)+(1-\alpha)\cos\gamma f'_{u}(r)=0, \quad \text{if} \quad u\neq -1/2.$$
			
			\item [(ii)] the radius $R_{sp}^{c}(g_{u})$ is the smallest positive root of the equation 
			$$rg''_{u}(r)+(1-\alpha)\cos\gamma g'_{u}(r)=0.$$
			
			\item [(iii)] the radius $R_{sp}^{c}(h_{u})$ is the smallest positive root of the equation 
			$$rh''_{u}(r)+(1-\alpha)\cos\gamma h'_{u}(r)=0.$$
		\end{enumerate}
	\end{theorem}
	\begin{proof}
		We begin with the first part. From \eqref{fL}, we have
		\begin{equation}\label{f1}
		1+\frac{zf''_{u}(z)}{f'_{u}(z)}= 1+\frac{z\mathcal{L}''_{u}(z)}{\mathcal{L}'_{u}(z)} +\left(\frac{1}{u+\frac{1}{2}}-1 \right)\frac{z\mathcal{L}'_{u}(z)}{\mathcal{L}_{u}(z)}.
		\end{equation}
		Also using the result \cite[Lemma~1]{Bricz-Rama}, we have
		\begin{equation*}
		\mathcal{L}_{u}(z)=\frac{z^{u+\tfrac{1}{2}}}{u(u+1)}\Phi_0(z)=\frac{z^{u+\frac{1}{2}}}{u(u+1)} \prod_{n\geq1}\left(1-\frac{z^2}{\tau^2_{u,n}}\right),
		\end{equation*}
		where $\Phi_k(z):={}_1F_{2}\left(1; \frac{u-k+2}{2}, \frac{u-k+3}{2};-\frac{z^2}{4}\right)$ with conditions as mentioned in \cite[Lemma~1]{Bricz-Rama}, and from the proof of \cite[Theorem~3]{Bricz-Rama}, we see that the entire function $\frac{u(u+1)}{u+\frac{1}{2}} z^{-u+\frac{1}{2}}\mathcal{L}'_{u}(z)$ is of order $1/2$ and therefore, has the following Hadamard factorization:
		\begin{equation*}
		\mathcal{L}'_{u}(z)=\frac{u+\frac{1}{2}}{u(u+1)}z^{u-\frac{1}{2}} \prod_{n\geq 1}\left(1-\frac{z^2}{\breve{\tau}^2_{u,n}} \right),
		\end{equation*}
		where $\tau_{u,n}$ and $\breve{\tau}_{u,n}$ are the $n$-th positive zeros of $\mathcal{L}_{u}$ and $\mathcal{L}'_{u}$, respectively and interlace for $0\neq u\in(-1,1)$ (see \cite[Theorem~1]{Bricz-Rama}). Now we can rewrite \eqref{f1} as follows:
		\begin{equation*}
		1+\frac{zf''_{u}(z)}{f'_{u}(z)}=1-\left(\frac{1}{u+\frac{1}{2}}-1\right)\sum_{n\geq 1}\frac{2z^2}{\tau^2_{u,n}-z^2}-\sum_{n\geq 1}\frac{2z^2}{\breve{\tau}^2_{u,n}-z^2}.
		\end{equation*}
		Let us now consider the case $u\in(0,1/2]$. Then using the inequality $||x|-|y||\leq |x-y|$ for $|z|=r<\breve{\tau}_{u,1}<\tau_{u,1}$ we get
		\begin{equation}\label{final-f}
		\left|\frac{zf''_{u}(z)}{f'_{u}(z)}\right|\leq \left(\frac{1}{u+\frac{1}{2}}-1\right)\sum_{n\geq 1}\frac{2r^2}{\tau^2_{u,n}-r^2}+\sum_{n\geq 1}\frac{2r^2}{\breve{\tau}^2_{u,n}-r^2} =-\frac{rf''_{u}(r)}{f'_{u}(r)}
		\end{equation}
		and for the case $u\in(1/2,1)$, using the inequality \eqref{firstnorm} with $\lambda=1-1/(u+1/2)$, we also get
		\begin{equation}\label{final-fu}
		\left|\frac{zf''_{u}(z)}{f'_{u}(z)}\right|\leq -\frac{rf''_{u}(r)}{f'_{u}(r)},
		\end{equation}	
		which is same as \eqref{final-f}. When $u\in (-1,0)$, then we proceed similarly substituting $u$ by $u-1$, $\Phi_0$ by $\Phi_1$, where $\Phi_1$ belongs to the Laguerre-P\'{o}lya class $\mathcal{LP}$ and the $n$-th positive zeros $\xi_{u,n}$ and $\breve{\xi}_{u,n}$ of $\Phi_1$ and its derivative $\Phi'_1$, respectively interlace. Finally, replacing $u$ by $u+1$, we obtain the required inequality.\\
		\indent For $0\neq u\in (-1,1)$, the Hadamard factorization for the entire functions $g'_{u}$ and $h'_{u}$ of order $1/2$ \cite[Theorem~3]{Bricz-Rama} is given by 
		\begin{equation}\label{Had-gh}
		g'_{u}(z)=\prod_{n\geq 1}\left( 1-\frac{z^2}{\gamma^2_{u,n}}\right) \quad \text{and} \quad
		h'_{u}(z)=\prod_{n\geq 1}\left( 1-\frac{z}{\delta^2_{u,n}}\right),
		\end{equation}
		where $\gamma_{u,n}$ and $\delta_{u,n}$ are $n$-th positive zeros of $g'_{u}$ and $h'_{u}$, respectively and $\gamma_{u,1}, \delta_{u,1} < \tau_{u,1}$. Now from \eqref{fL} and \eqref{Had-gh}, we have
		\begin{equation}\label{gh-u}
		\left\{
		\begin{array}{lr}
		1+\frac{zg''_{u}(z)}{g'_{u}(z)}= \frac{1}{2}-u +z\dfrac{(\frac{3}{2}-u)\mathcal{L}'_{u}(z)+z\mathcal{L}''_{u}(z)}{(\frac{1}{2}-u)\mathcal{L}_{u}(z)+z\mathcal{L}'_{u}(z)} =1-\sum_{n\geq 1}\dfrac{2z^2}{\gamma^2_{u,n}-z^2} \\
		1+\frac{zh''_{u}(z)}{h'_{u}(z)}= \frac{1}{2}\left(\frac{3}{2}-u +\sqrt{z}\dfrac{(\frac{5}{2}-u)\mathcal{L}'_{u}(\sqrt{z})+\sqrt{z}\mathcal{L}''_{u}(\sqrt{z})}{(\frac{3}{2}-u)\mathcal{L}_{u}(\sqrt{z})+\sqrt{z}\mathcal{L}'_{u}(\sqrt{z})}    \right)=1-\sum_{n\geq 1}\dfrac{z}{\delta^2_{u,n}-z}.
		\end{array}
		\right.
		\end{equation}
		Using the inequality $||x|-|y||\leq |x-y|$ in \eqref{gh-u} for $|z|=r< \gamma_{u,1}$ and $|z|=r<\delta_{u,1}$, we get
		\begin{equation}\label{final-gh}
		\left\{
		\begin{array}{lr}
		\left|\frac{zg''_{u}(z)}{g'_{u}(z)} \right|\leq \sum_{n\geq 1}\frac{2r^2}{\gamma^2_{u,n}-r^2}=-\frac{rg''_{u}(r)}{g'_{u}(r)} \\
		\left|\frac{zh''_{u}(z)}{h'_{u}(z)} \right|\leq \sum_{n\geq 1}^{\infty}\frac{r}{\delta^2_{u,n}-r}=-\frac{rh''_{u}(r)}{h'_{u}(r)}. 
		\end{array}
		\right.
		\end{equation} 
		Further, proceeding with the similar method as in Theorem~\ref{wright-phi}, result follows. \qed
	\end{proof}	
	 
	 With similar reasoning as Theorem~\ref{wright-phi}, the proof of the following holds.
	\begin{theorem}\label{lommel-convexgamma}
		Let $u\in (-1,1)$, $u\neq0$. Then the radius of $\gamma$-Spirallikeness of order $\alpha$
		for the functions Let the functions $f_{u}, g_{u}$ and $h_{u}$ given by \eqref{lomel-normalized} are the smallest positive roots of the following equations:
		\begin{enumerate}
			\item [$(i)$] $rf'_{u}(r)+((1-\alpha)\cos\gamma-1)f_{u}(r)=0$
			\item [$(ii)$] $rg'_{u}(r)+((1-\alpha)\cos\gamma-1)g_{u}(r)=0$
			\item [$(iii)$] $rh'_{u}(r)+((1-\alpha)\cos\gamma-1)h_{u}(r)=0$
		\end{enumerate} 
	in $(0,\tau_{u,1})$, $(0,\tau_{u,1})$ and $(0,\tau^2_{u,1})$, respectiively.
	\end{theorem}
	\begin{remark}
		Taking $\gamma=0$, Theorem~\ref{lommel-convexgamma} reduces to \cite[Theorem~3]{Bricz-Rama}.
	\end{remark}	
	
	\section{Struve functions}
	The Struve function $\mathcal{\bf{H}}_{\beta}$ of first kind is a particular solution of the second-order inhomogeneous Bessel differential equation $$z^2w''(z)+zw'(z)+(z^2-{\beta}^2)w(z)=\frac{4\left(\frac{z}{2}\right)^{\beta+1}}{\sqrt{\pi}\Gamma\left(\beta+\frac{1}{2}\right)}$$ and have the following form:
	\begin{equation*}
	\mathcal{\bf{H}}_{\beta}(z):=\frac{\left(\frac{z}{2}\right)^{\beta+1}}{\sqrt{\frac{\pi}{4}}\Gamma\left(\beta+\frac{1}{2}\right)} {}_1 F_{2}\left(1;\frac{3}{2},\beta+\frac{3}{2};-\frac{z^2}{4}\right) ,
	\end{equation*}
	where $-\beta-\frac{3}{2}\notin\mathbb{N}$ and ${}_1 F_{2}$ is a hypergeometric function. Since it is not normalized, therefore we take the normalized functions:
	\begin{equation}\label{nor-strv-uvw}
	\left\{
	\begin{array}{lr}
	U_{\beta}(z)=\left(\sqrt{\pi}2^{\beta}\left(\beta+\frac{3}{2}\right){\bf{H}}_{\beta}(z)\right)^{\frac{1}{\beta+1}},\\
	V_{\beta}(z)=\sqrt{\pi}2^{\beta}z^{-\beta}\Gamma\left(\beta+\frac{3}{2}\right){\bf{H}}_{\beta}(z),\\
	W_{\beta}(z)=\sqrt{\pi}2^{\beta}z^{\frac{1-\beta}{2}}\Gamma\left(\beta+\frac{3}{2}\right){\bf{H}}_{\beta}(\sqrt{z}).
	\end{array}
	\right.
	\end{equation}
	Moreover, for $|\beta|\leq{1}/{2}$, ${\bf{H}}_{\beta}$ (see \cite[Lemma~1]{bricz-2017}) and ${\bf{H}}'_{\beta}$ have the Hadamard factorizations \cite[Theorem~4]{Bricz-Rama} given by
	\begin{equation*}
	{\bf{H}}_{\beta}(z)=\frac{z^{\beta+1}}{\sqrt{\pi}2^{\beta}\Gamma(\beta+\frac{3}{2})}\prod_{n\geq1}\left(1-\frac{z^2}{z^2_{\beta,n}}\right)
	\end{equation*}
	and
	\begin{equation}\label{strv-facto}
	{\bf{H}}'_{\beta}(z)=\frac{(\beta+1)z^{\beta}}{\sqrt{\pi}2^{\beta}\Gamma(\beta+\frac{3}{2})}\prod_{n\geq1}\left(1-\frac{z^2}{\breve{z}^2_{\beta,n}}\right)
	\end{equation}
	where $z_{\beta,n}$ and $\breve{z}_{\beta,n}$ are the $n$-th positive zeros of ${\bf{H}}_{\beta}$ and ${\bf{H}}'_{\beta}$,respectively and interlace \cite[Theorem~2]{Bricz-Rama}. Thus from \eqref{strv-facto} with logarithmic differentiation, we obtain respectively
	\begin{equation*}
	\frac{z{\bf{H}}'_{\beta}(z)}{{\bf{H}}_{\beta}(z)}=(\beta+1)-\sum_{n\geq1}\frac{2z^2}{z^2_{\beta,n}-z^2} 
	\end{equation*}
	and
	\begin{equation}\label{strv-strlikeconvex}
	1+\frac{z{\bf{H}}''_{\beta}(z)}{{\bf{H}}'_{\beta}(z)}=(\beta+1)-\sum_{n\geq1}\frac{2z^2}{\breve{z}^2_{\beta,n}-z^2} .
	\end{equation}
	Also for $|\beta|\leq {1}/{2}$, the Hadamard factorization for the entire functions $V'_{\beta}$ and $W'_{\beta}$ of order $1/2$ \cite[Theorem~4]{Bricz-Rama} is given by 
	\begin{equation}\label{Had-VW'}
	V'_{\beta}(z)=\prod_{n\geq 1}\left( 1-\frac{z^2}{\eta^2_{\beta,n}}\right) \quad \text{and} \quad
	W'_{\beta}(z)=\prod_{n\geq 1}\left( 1-\frac{z}{\sigma^2_{\beta,n}}\right),
	\end{equation}
	where $\eta_{\beta,n}$ and $\sigma_{\beta,n}$ are $n$-th positive zeros of $V'_{\beta}$ and $W'_{\beta}$, respectively. $V'_{\beta}$ and $W'_{\beta}$ belong to the Laguerre-P\'{o}lya class and zeros satisfy $\eta_{\beta,1}, \sigma_{\beta,1} < z_{\beta,1}$. Now proceeding as in Theorem~\ref{lommeltheorem} using \eqref{nor-strv-uvw}, \eqref{strv-facto}, \eqref{strv-strlikeconvex} and \eqref{Had-VW'}, we obtain the following results:

	\begin{theorem}\label{struve-gamma}
	Let $|\beta|\leq 1/2$. Then the radii of $\gamma$-Spirallikeness of order $\alpha$
	for the functions $U_{\beta}, V_{\beta}$ and $W_{\beta}$ given by \eqref{nor-strv-uvw} are the smallest positive roots of the following equations:
	\begin{enumerate}
		\item [(i)] $rU'_{\beta}(r)+((1-\alpha)\cos\gamma-1)U_{\beta}(r)=0$
		\item [(ii)] $rV'_{\beta}(r)+((1-\alpha)\cos\gamma-1)V_{\beta}(r)=0$
		\item [(iii)] $rW'_{\beta}(r)+((1-\alpha)\cos\gamma-1)W_{\beta}(r)=0$
	\end{enumerate}
	in $(0, z_{\beta,1})$, $(0,z_{\beta,1})$ and $(0,z^2_{\beta,1})$, respectively.
\end{theorem}
\begin{remark}
		Taking $\gamma=0$ in Theorem~\ref{struve-gamma} gives \cite[Theorem~2]{bdoy-2016}.
\end{remark}
		
	\begin{theorem}\label{struvetheorem}
		Let $|\beta|\leq 1/2$. Let the functions $U_{\beta}, V_{\beta}$ and $W_{\beta}$ be given by \eqref{nor-strv-uvw}. Then
		\begin{enumerate}
			\item [(i)] the radius $R_{sp}^{c}(U_{\beta})$ is the smallest positive root of the equation 
			$$rU''_{\beta}(r)+(1-\alpha)\cos\gamma U'_{\beta}(r)=0.$$
			
			\item [(ii)] the radius $R_{sp}^{c}(V_{\beta})$ is the smallest positive root of the equation 
			$$rV''_{\beta}(r)+(1-\alpha)\cos\gamma V'_{\beta}(r)=0.$$
			
			\item [(iii)] the radius $R_{sp}^{c}(W_{\beta})$ is the smallest positive root of the equation 
			$$rW''_{\beta}(r)+(1-\alpha)\cos\gamma W'_{\beta}(r)=0.$$
		\end{enumerate}
	\end{theorem}
\begin{remark}
	Taking $\gamma=0$ in Theorem~\ref{struvetheorem} gives \cite[Theorem~4]{Bricz-Rama}.
\end{remark}

	\section{On Ramanujan type entire functions}
	Ismail and Zhang~\cite{ismail-2018} defined the following entire function of growth order zero for $\beta>0$, called Ramanujan type entire function
	\begin{equation*}
	A^{(\beta)}_{p}(c,z)= \sum_{n\geq0}\frac{(c;p)_n p^{\beta n^2}}{(p;p)_n}z^n,
	\end{equation*}
	where $\beta>0$, $0<p<1$, $c\in \mathbb{C}$, $(c;p)_0=1$ and $(c;p)_k=\prod_{j=0}^{k-1}(1-cp^j)$ for $k\geq1,$ which is the generalization of both the Ramanujan entire function $A_{p}(z)$ and Stieltjes-Wigert polynomial $S_n(z;p)$ defined as (see \cite{ismail-2005,ramanujan-1988}):
	\begin{equation*}
	A_{p}(-z)=A^{(1)}_{p}(0,z)= \sum_{n=0}^{\infty}\frac{p^{n^2}}{(p;p)_n}z^n  
	\end{equation*}
	and
	\begin{equation*}
	A^{(1/2)}_{p}(p^{-n},z)=\sum_{m=0}^{\infty}\frac{(p^{-n};p)_m p^{m^2/2}}{(p;p)_m} z^m = (p;p)_n S_n(zp^{(1/2)-n};p). 
	\end{equation*}
	Since $A^{(\beta)}_{p}(c,z)\not \in \mathcal{A}$, therefore consider the following three normalized functions in $\mathcal{A}$:
	\begin{equation}\label{ramj1}
	\left\{
	\begin{array}{lr}
	f_{\beta, p,c}(z):= \left(z^{\beta} A^{(\beta)}_{p}(-c,-z^2)  \right)^{1/\beta} \\
	g_{\beta, p,c}(z):= z A^{(\beta)}_{p}(-c,-z^2) \\
	h_{\beta, p,c}(z):= z A^{(\beta)}_{p}(-c,-z),
	\end{array}
	\right.
	\end{equation}
	where $\beta>0$, $c\geq0$ and $0<p<1$. From \cite[Lemma~2.1, p.~4-5]{ErhanDenij2020}, we see that the function
	$$z\rightarrow \Psi_{\beta,p,c}(z):= A^{(\beta)}_{p}(-c,-z^2)$$
	has infinitely many zeros (all are positive) for $\beta>0$, $c\geq0$ and $0<p<1$. Let $\psi_{\beta,p,n}(c)$ be the $n$-th positive zero of $\Psi_{\beta,p,c}(z)$. Then it has the following Weiersstrass decomposition:
	\begin{equation}\label{ramj2}
	\Psi_{\beta,p,c}(z)= \prod_{n\geq1}\left(1-\frac{z^2}{\psi^2_{\beta, p, n}(c)} \right).
	\end{equation}
	Moreover, the $n$-th positive zero $\Xi_{\beta, p, n}(c)$ of the derivative of the  following function
	\begin{equation}\label{ramj3}
	\Phi_{\beta, p,c}(z):= z^{\beta}\Psi_{\beta, p,c}(z)
	\end{equation}
	interlace with $\psi_{\beta, p, n}(c)$ and satisfy the relation
	$$\Xi_{\beta, p, n}(c)< \psi_{\beta,p,n}(c)<\Xi_{\beta, p, n+1}(c)< \psi_{\beta,p,n+1}(c)$$ 
	for $n\geq1.$ Now using \eqref{ramj1} and \eqref{ramj2}, we have
	\begin{align*}\label{ramj-star}
	\frac{zf'_{\beta,p,c}(z)}{f_{\beta,p,c}(z)}&= 1+ \frac{1}{\beta}\frac{z\Psi'_{\beta,p,c}(z)}{\Psi_{\beta,p,c}(z)}
	=1-\frac{1}{\beta}\sum_{n\geq1}\frac{2z^2}{\psi^2_{\beta,p,n}(c)- z^2} ; \;(c>0) \nonumber\\
	\frac{zg'_{\beta,p,c}(z)}{g_{\beta,p,c}(z)}&= 1+ \frac{z\Psi'_{\beta,p,c}(z)}{\Psi_{\beta,p,c}(z)}
	=1-\sum_{n\geq1}\frac{2z^2}{\psi^2_{\beta,p,n}(c)- z^2} ; \nonumber\\
	\frac{zh'_{\beta,p,c}(z)}{h_{\beta,p,c}(z)}&= 1+ \frac{1}{2}\frac{\sqrt{z} \Psi'_{\beta,p,c}(\sqrt{z})}{\Psi_{\beta,p,c}(\sqrt{z})}
	=1-\sum_{n\geq1}\frac{z}{\psi^2_{\beta,p,n}(c)- z},
	\end{align*}
	where $\beta>0,c\geq0$ and $0<p<1$. Also, using \eqref{ramj3} and the infinite product representation of $\Phi'$~\cite[p.~14-15, Also see Eq.~4.6]{ErhanDenij2020}, we have
	\begin{align*}
	1+\frac{zf''_{\beta,p,c}(z)}{f'_{\beta,p,c}(z)}&= 1+\frac{z \Phi''_{\beta,p,c}(z)}{\Phi'_{\beta,p,c}(z)} 
	+\left(\frac{1}{\beta}-1\right) \frac{z \Phi'_{\beta,p,c}(z)}{\Phi_{\beta,p,c}(z)}\\
	&=1-\sum_{n\geq1}\frac{2z^2}{\Xi^2_{\beta,p,n}(c)- z^2} -\left(\frac{1}{\beta}-1\right) \sum_{n\geq1}\frac{2z^2}{\psi^2_{\beta,p,n}(c)- z^2}.
	\end{align*}
	As $(z\Psi_{\beta,p,c}(z))'$ and $h'_{\beta,p,c}(z)$ belong to $\mathcal{LP}$. So suppose $\gamma_{\beta,p,n}(c)$ be the positive zeros of $g'_{\beta,p,c}(z)$ (growth order is same as $\Psi_{\beta,p,c}(z)$) and $\delta_{\beta,p,n}(c)$ be the positive zeros of $h'_{\beta,p,c}(z)$. Thus using their infinite product representations, we have
	\begin{align*}
	1+\frac{zg''_{\beta,p,c}(z)}{g'_{\beta,p,c}(z)}&= 1-\sum_{n\geq1}\frac{2z^2}{\gamma^2_{\beta,p,n}(c)- z^2} \nonumber\\
	1+\frac{zh''_{\beta,p,c}(z)}{h'_{\beta,p,c}(z)}&= 1-\sum_{n\geq1}\frac{z}{\delta^2_{\beta,p,n}(c)- z}. 
	\end{align*}
	Now proceeding similarly as done in the above sections, we obtain the following results:
	\begin{theorem}\label{ramjThm1}
		Let $\beta>0$, $c\geq0$ and $0<p<1$. Then the radii of $\gamma$-Spirallikeness of order $\alpha$ for the functions $f_{\beta, p,c}(z)$, $g_{\beta, p,c}(z)$ and $h_{\beta, p,c}(z)$ given by \eqref{ramj1} are the smallest positive roots of the following equations:
		\begin{enumerate}
			\item [(i)] $rf'_{\beta, p,c}(r)+((1-\alpha)\cos\gamma-1)f_{\beta, p,c}(r)=0$
			\item [(ii)] $rg'_{\beta, p,c}(r)+((1-\alpha)\cos\gamma-1)g_{\beta, p,c}(r)=0$
			\item [(iii)] $rh'_{\beta, p,c}(r)+((1-\alpha)\cos\gamma-1)h_{\beta, p,c}(r)=0$
		\end{enumerate}
	in $(0, \psi_{\beta,p,1}(c))$, $(0,\psi_{\beta,p,1}(c))$ and $(0,\psi^2_{\beta,p,1}(c))$, respectively.
	\end{theorem}
	
	 We now conclude this section with the convex analog of Theorem~\ref{ramjThm1}.
	\begin{theorem}\label{ramjThm2}
		Let $\beta>0$, $c\geq0$ and $0<p<1$. Let the functions $f_{\beta, p,c}(z)$, $g_{\beta, p,c}(z)$ and $h_{\beta, p,c}(z)$ be given by \eqref{ramj1}. Then
		\begin{enumerate}
			\item [(i)] the radius $R_{sp}^{c}(f_{\beta, p,c}(z))$ is the smallest positive root of the equation 
			$$rf''_{\beta, p,c}(r)+(1-\alpha)\cos\gamma f'_{\beta, p,c}(r)=0.$$
			
			\item [(ii)] the radius $R_{sp}^{c}(g_{\beta, p,c}(z))$ is the smallest positive root of the equation 
			$$rg''_{\beta, p,c}(r)+(1-\alpha)\cos\gamma g'_{\beta, p,c}(r)=0.$$
			
			\item [(iii)] the radius $R_{sp}^{c}(h_{\beta, p,c}(z))$ is the smallest positive root of the equation 
			$$rh''_{\beta, p,c}(r)+(1-\alpha)\cos\gamma h'_{\beta, p,c}(r)=0.$$
		\end{enumerate}
	\end{theorem}

	%\begin{acknowledgements}
	%We are thankful to the Editor and the Reviewers for their valuable suggestions to improve the previous version of this manuscript. 
	%\end{acknowledgements}

	% Authors must disclose all relationships or interests that 
	% could have direct or potential influence or impart bias on 
	% the work: 
	%
	\section*{Statements and Declarations}
	\begin{itemize}
		%	\item {\bf Funding}: The work of Kamajeet Gangania is supported by University Grant Commission, New-Delhi, India  under UGC-Ref. No.:1051/(CSIR-UGC NET JUNE 2017).
		\item {\bf Conflict of interest}: The authors declare that they have no conflict of interest
		%	\item Ethics approval 
		%	\item Consent to participate
		%	\item Consent for publication
		\item {\bf Availability of data and materials }: None
		%	\item Code availability 
		\item {\bf Authors' contributions }: All authors contributed Equally.
	\end{itemize}

	% BibTeX users please use one of
	%\bibliographystyle{spbasic}      % basic style, author-year citations
	%\bibliographystyle{spmpsci}      % mathematics and physical sciences
	%\bibliographystyle{spphys}       % APS-like style for physics
	%\bibliography{}   % name your BibTeX data base

\begin{thebibliography}{20}
		%
		% and use \bibitem to create references. Consult the Instructions
		% for authors for reference list style.
		%
		%\bibitem{RefJ}
		%% Format for Journal Reference
		%Author, Article title, Journal, Volume, page numbers (year)
		%% Format for books
		%\bibitem{RefB}
		%Author, Book title, page numbers. Publisher, place (year)
		%% etc
		
		\bibitem{abo-2018}  Akta\c{s},  \.{I}. Baricz, \'{A}. and Orhan, H.: Bounds for radii of starlikeness and convexity of some special functions. {Turkish J. Math. } {\bf 42}, 211--226 (2018) doi: 10.3906/mat-1610-41
		
		%\bibitem{baric-q2016} { Baricz, \'{A}. Dimitrov, D.K. and Mez\H{o}, I.}: Radii of starlikeness and convexity of some $q$-Bessel functions. {J. Math. Anal. Appl.} {\bf 435}, 968--985 (2016). doi: 10.1016/j.jmaa.2015.10.065
		
		\bibitem{bdoy-2016}  {Baricz, \'{A}. Dimitrov, D.K., Orhan, H. and Ya\u{g}mur, N.}: Radii of starlikeness of some special functions. {Proc. Amer. Math. Soc.} { \bf 144}, 3355--3367 (2016). doi: 10.1090/proc/13120
		
		\bibitem{Baric-2014} { Baricz, \'{A}. Kup\'{a}n, P.A. and Sz\'{a}sz, R.}: The radius of starlikeness of normalized Bessel functions of the first kind. {Proc. Amer. Math. Soc.} {\bf 142}, 2019--2025 (2014). doi: 10.1090/S0002-9939-2014-11902-2
		
		%\bibitem{baricz-2010} { Baricz, \'{A}. and Ponnusamy, S.}: Starlikeness and convexity of generalized Bessel functions. {Integral Transforms Spec. Funct. } {\bf 21}, 641--653 (2010). doi: 10.1080/10652460903516736
		
		\bibitem{bricz-2017} { Baricz, \'{A}. Ponnusamy, S. and Singh, S.}: Tur\'{a}n type inequalities for Struve functions. {J. Math. Anal. Appl.} {\bf 445}, 971--984 (2017). doi: 10.1016/j.jmaa.2016.08.026
		
		%\bibitem{Bansel-2016} { Bansal, D. and Prajapat, J.K.}: Certain geometric properties of the Mittag-Leffler functions. { Complex Var. Elliptic Equ.} {\bf 61}, 338--350 (2016). doi: 10.1080/17476933.2015.1079628
		
		\bibitem{b-praj-2020} { Baricz, \'{A}. and Prajapati, A.}: Radii of starlikeness and convexity of generalized Mittag-Leffler functions. {Math. Commun.} {\bf 25}, 117--135 (2020).
		
		\bibitem{Baric-2015} { Baricz, \'{A}. and Sz\'{a}sz, R.}: The radius of convexity of normalized Bessel functions. {Anal. Math.} { \bf 41}, 141--151 (2015). doi: 10.1007/s10476-015-0202-6
		
		\bibitem{btk-2018} { Baricz, \'{A}. Toklu, E. and Kadio\u{g}lu, E.}: Radii of starlikeness and convexity of Wright functions. {Math. Commun.} {\bf 23}, 97--117 (2018).
		
		\bibitem{Bricz-Rama} { Baricz, \'{A}. and Ya\u{g}mur, N.}: Geometric properties of some Lommel and Struve functions. {Ramanujan J. } {\bf 42}, 325--346 (2017). doi: 10.1007/s11139-015-9724-6
		
		\bibitem{bulut-engel-2019} {Bulut, S. and Engel, O.}: The radius of starlikeness, convexity and uniform convexity of the Legendre polynomials of odd degree. {Results Math.} { \bf74}, Paper No. 48, 9 pp (2019). doi: 10.1007/s00025-019-0975-1
		
		\bibitem{ErhanDenij2020} { Deniz, E.}: Geometric and monotonic properties of Ramanujan type entire functions. {Ramanujan J.} {\bf 55}, 103--130 (2020). doi: 10.1007/s11139-020-00267-w
		
		\bibitem{Deniz-2017} { Deniz, E. and  Sz\'{a}sz, R.}: The radius of uniform convexity of Bessel functions. {J. Math. Anal. Appl.} {\bf453}, 572--588 (2017). doi: 10.1016/j.jmaa.2017.03.079
		
		\bibitem{lp} { Dimitrov, D.K. and Ben Cheikh, Y.}. Laguerre polynomials as Jensen polynomials of Laguerre-P\'{o}lya entire functions. {J. Comput. Appl. Math.} { \bf 233}, 703--707 (2009). doi: 10.1016/j.cam.2009.02.039
		
		\bibitem{g-specialIJST} Gangania, K. and Kumar, S.S.: $\mathcal{S}^*(\psi)$ and $\mathcal{C}(\psi)$-radii for some special functions. Iran. J. Sci. Technol. Trans. A Sci. {\bf 46}, 955–966 (2022).
		
		%\bibitem{ganga1997} { Gangadharan, A., Ravichandran, V. and Shanmugam, T.N.}: Radii of convexity and strong starlikeness for some classes of analytic functions. {J. Math. Anal. Appl. } {\bf 211} 301--313 (1997). doi: 10.1006/jmaa.1997.5463
		
		
		\bibitem{ismail-2018} { Ismail, M.E.H. and Zhang, R.}: $q$-Bessel functions and Rogers-Ramanujan type identities. {Proc. Amer. Math. Soc.} {\bf 146}, 3633--3646 (2018). doi: 10.1090/proc/13078
		
		\bibitem{ismail-2005} { Ismail, M.E.H.}: Classical and quantum orthogonal polynomials in one variable. {\em Encyclopedia of Mathematics and its Applications}, 98, Cambridge University Press, Cambridge, 2005.
		
		\bibitem{SG-2020}  Kumar, S.S. and Gangania, K.: Subordination and radius problems for certain starlike functions. arXiv:2007.07816 (2020)
					
		\bibitem{pathan-2016} { Kumar H and Pathan A.M}. On the distribution of non-zero zeros of generalized Mittag-Leffler functions. {International Journal of Engineering Research and Application}, { \bf 6}, 66--71 (2016).
		
		\bibitem{Levin-1996}  { Ya. Levin B}. Lectures on Entire Functions. { Translation of Mathematics Monographs, American Mathematical Society}, Providence 150 (1996).
		
		\bibitem{minda94}  { Ma, W.C. and Minda, D.}: A unified treatment of some special classes of univalent functions. { Proceedings of the Conference on Complex Analysis, Tianjin}, Conf Proc Lecture Notes Anal, I Int Press. Cambridge, MA. 157--169 (1992). 
				
        \bibitem{Pfaltzgraff-1975} Pfaltzgraff, J.A.:  Univalence of the integral of $f\sp{\prime} (z)\sp{\lambda }$. Bull. London Math. Soc. {\bf 7}, no.~3, 254--256  (1975).
        	
		\bibitem{prabha-1971} { Prabhakar, T.R.}: A singular integral equation with a generalized Mittag Leffler function in the kernel. {Yokohama Math. J.} {\bf 19}: 7--15 (1971).
		
		\bibitem{ramanujan-1988} { Ramanujan, S.}: The lost notebook and other unpublished papers, Springer-Verlag, Berlin, 1988.
		
		\bibitem{Robertson-1969}  Robertson, M.S.: Univalent functions $f(z)$ for which $zf\sp{\prime} (z)$ is spirallike. Michigan Math. J. {\bf 16}, 97--101 (1969).
		
		\bibitem{spacek-1933} Spacek, L.: Contribution \'{a} la th\`{e}orie des fonctions univalentes, Casop Pest. Mat.-Fys. {\bf 62}, 12-19  (1933).
		 
		\bibitem{Kazimoglu-2022} Kazımoğlu, S. and  Deniz, E.: Radius Problems for Functions Containing Derivatives of Bessel Functions. Comput. Methods Funct. Theory (2022). https://doi.org/10.1007/s40315-022-00455-3 

				
	%	\bibitem{Toklu-q2019} { Toklu, E.}: Radii of starlikeness and convexity of $q$-Mittag-Leffler functions. {Turkish J. Math.} {\bf 43}, 2610--2630 (2019). doi: 10.3906/mat-1907-54
		
		\bibitem{watson-1944} { Watson, G.N. }: A Treatise on the Theory of Bessel Functions, Cambridge University Press, 1944.
		
		
		
		
	\end{thebibliography}
	%
	% Non-BibTeX users please use

\end{document}